\documentclass[12pt]{article}
\usepackage{amssymb,amsmath}
\newtheorem{lemma}{Lemma}[section]
\newtheorem{theorem}{Theorem}[section]

\newtheorem{corollary}{Corollary}[section]

\newcounter{saveeqn}%
\begin{document}
\title{ Sharp  Nash Inequalities on the unit sphere\\
         The influence of symmetries}
\author{{\Large Athanase Cotsiolis and Nikos Labropoulos}\\
Department of Mathematics, University of Patras,\\ Patras 26110, Greece\\
\small{e-mails: cotsioli@math.upatras.gr and nal@upatras.gr}}
\maketitle

{\noindent}{\bf Abstract:}  {{  In this paper both we establish
the best constants for the Nash inequalities on the standard unit
sphere $\mathbb{S}^n$ of $\mathbb{R}^{n+1}$ and we give answers on
the existence of extremal functions on the corresponding problems.
Also we study the problem of the best constants in the case, where
the data are invariant under the
action of the group $G=O(k)\times O(m)$, and we find the best constants.\\
{\noindent}{\textbf{\emph{Keywords:}} Manifolds without boundary;
Standard unit sphere; Nash inequalities; Best constants; Extremal
functions; Symmetries.}

\section{Introduction}
\setcounter{equation}{0}$\;\;\;\;$Nash inequalities after their
first appearance in the celebrated paper of Nash \cite{Nas},
reappear in some subsequent papers. Specifically, we refer to
\cite{Car-Los, Bec, Dru-Heb-Vau, Hum1} for manifolds without
boundary and \cite{Hum2, Hum3, Cot-Lab-Tra, Cot-Lab} for manifolds
with boundary. In this paper we are focusing our interest in the
special case when the manifold is the standard unit sphere
$\mathbb{S}^n$ of $\mathbb{R}^{n+1}$.

Let $(M,g)$ be a smooth, complete $n-$dimensional Riemannian
manifold of infinite volume, where $n\geq 1$.\\
We say that the Nash inequality (\ref{E1.1}) is valid if there
exists a constant $A>0$ such that for all $u \in C_0^\infty \left(
M \right)$,
\begin{equation}\label{E1.1}
 \left( {\int_M {u^2 dV_g } }
\right)^{1 + \frac{2} {n}} \leqslant A  {\int_M {\left| {\nabla u}
\right|_g^2 } dV_g   } \left( {\int_M {\left| u \right|dV_g } }
\right)^{\frac{4} {n}}
\end{equation}
Such an inequality, as refereed above,  first appeared in the
celebrated paper of Nash \cite{Nas}, when discussing the H\"older
regularity of solutions of divergence form
uniformly elliptic equations.\\
Let $A_0(n)$ be the best constant in Nash's inequality
(\ref{E1.1}) above for the Euclidean space. That is
\[
A_0 \left( n \right)^{ - 1}  = \inf \left\{ {\left.
{\frac{{\int_{\mathbb{R}^n } {\left| {\nabla u} \right|^2 dx}
\left( {\int_{\mathbb{R}^n } {\left| u \right|dx} }
\right)^{\frac{4} {n}} }} {{\left( {\int_{\mathbb{R}^n } {u^2 dx}
} \right)^{1 + \frac{2} {n}} }}} \right|u \in C_0^\infty  \left(
{\mathbb{R}^n } \right),u\not \equiv 0} \right\}
\]
This best constant has been computed by Carlen and Loss in
\cite{Car-Los}, together with the characterization of the
extremals for the corresponding optimal inequality, as
\[
A_0 \left( n \right) = \frac{{\left( {n + 2} \right)^{\frac{{n +
2}} {n}} }} {{2^{\frac{2} {n}} n\lambda _{1,n} \left|
{\mathcal{B}^n } \right|^{\frac{2} {n}} }},
\]
where $\left| {\mathcal{B}^n } \right|$ denotes the euclidian
volume of the unit ball $\mathcal{B}^n $ in $\mathbb{R}^n $ and
$\lambda _{1,n} $ is the first Neumann eigenvalue for the
Laplacian for radial functions in the  unit ball $\mathcal{B}^n
$.\\
For an example of application of the Nash inequality with the best
constant, we refer to Kato \cite{Kat} and for a geometric proof
with an asymptotically sharp constant, we refer to Beckner
\cite{Bec}.\\

For compact Riemannian manifolds, or smooth bounded domains, (see
Nirenberg \cite{Nir}), the Nash inequality still holds with an
additional $L^1-$term and that is why we will refer to this as the
$L^1-$Nash inequality.\\
Given $(M,g)$ a smooth compact $n-$dimensional Riemannian
manifold, $n\geq 2$, we are looking for the existence of real
constants $A$ and $B$ such that for any $u \in C^\infty(M)$,
\begin{equation}\label{E1.2}
 \left( {\int_M {u^2 dV_g } }
\right)^{1 + \frac{2} {n}} \leqslant A  {\int_M {\left| {\nabla u}
\right|_g^2 } dV_g   } \left( {\int_M {\left| u \right|dV_g } }
\right)^{\frac{4} {n}}+B\left( {\int_M {\left| u \right|dV_g } }
\right)^{2+\frac{4} {n}}
\end{equation}
One can define
\[
A^1_{opt}(M)  = \inf \left\{ {A > 0:\exists \,B >
0\,\,\,\mathrm{s.t.}\,\,\,(\ref{E1.2})
\,\,\,\mathrm{is\,\,true}\,\,\,\forall \,u \in C^\infty  \left( M
\right)} \right\}
 \]
 and
\[
B^1_{opt}(M)  = \inf \left\{ {B > 0:\exists \,A >
0\,\,\,\mathrm{s.t.}\,\,\,(\ref{E1.2})
\,\,\,\mathrm{is\,\,true}\,\,\,\forall \,u \in C^\infty  \left( M
\right)} \right\}
 \]
Druet, Hebey and Vaugon proved in \cite {Dru-Heb-Vau} that
$A^1_{opt}(M)=A_0(n)$, and (\ref{E1.2}) with its optimal constant
$A=A_0(n)$ is sometimes valid and sometimes not, depending on the
geometry, specifically on the sign of the curvature. This is
another illustration of the important idea of Druet \cite{Dru}
that an inequality may be at the same time localisable and
affected by the geometry. On the contrary,
$B^1_{opt}(M)=Vol(M)^{-1-2/n}$, where $Vol(M)$ is the volume of
the manifold, and (\ref{E1.2}) with its optimal constant
$B^1_{opt}(M)=Vol(M)^{-1-2/n}$ is always valid with geometry
playing no role (see also \cite {Dru-Heb-Vau}).\\
For all $u \in C^\infty  \left( M \right)$, consider now the
$L^2-$Nash inequality
\begin{equation}\label{E1.3}
 \left( {\int_M {u^2 dV_g } }
\right)^{1 + \frac{2} {n}} \leqslant \left( {A\int_M {\left|
{\nabla u} \right|_g^2 } dV_g  + B\int_M {u^2 dV_g } }
\right)\left( {\int_M {\left| u \right|dV_g } } \right)^{\frac{4}
{n}}
\end{equation}
 and define
\[
A^2_{opt}(M)  = \inf \left\{ {A > 0:\exists \,B >
0\,\,\,\mathrm{s.t.}\,\,\,(\ref{E1.3})
\,\,\,\mathrm{is\,\,true}\,\forall \,u \in C^\infty  \left( M
\right)} \right\}
 \]
 and
$$
B^2_{opt}(M)  = \inf \left\{ {B > 0:\exists \,A >
0\,\,\,\mathrm{s.t.}\,\,\,(\ref{E1.3})
\,\,\,\mathrm{is\,\,true}\,\,\,\forall \,u \in C^\infty  \left( M
\right)} \right\}
$$
Humbert studied in \cite{Hum1}  the $L^2-$Nash inequality in
detail. Contrary to the sharp $L^1-$Nash inequality, he proved in
this case that $B$ always exists and $A^2_{opt}(M)=A_0(n)$. Also,
he studied the second optimal constant $B^2_{opt}(M)$ of this
inequality, giving its explicit value
$B^2_{opt}(\mathbb{S}^1)=(2\pi)^{-2}$ for $n=1$ (i.e. for
$M=\mathbb{S}^1$), and, for $n>1$, proving that
$$
B^2_{opt}\geq\max\!\left(Vol(M)^{-2/n},\frac{|\mathcal{B}|^{-2/n}}
{6n}\left(\frac{2}{n+2}+\frac{n-2}{\lambda_1}\right)\!
\left(\frac{n+2}{2}\right)^{2/n}\! \!\max_{x\in M}S_g(x)\!
\right),
$$
where $|\mathcal{B}|$ is the volume of the unit ball $\mathcal{B}$
in $\mathbb{R}^n$, $\lambda_1$ is the first non-zero Neumann
eigenvalue of the Laplacian on radial functions on $\mathcal{B}$,
$Vol(M)$ is the volume of $(M,g)$ and $S_g(x)$ is the scalar
curvature of $g$ at $x$. In the same paper it was proved that, if
$(M,g)$ is a smooth compact Riemaniann $n-$manifold with $n\geq 1$
and $L_1-$Nash inequality is true, with $A=A_0(n)$ and some $B$,
then there exists $u_0\in H_1(M),\,\, u_0\not  \equiv 0$, (where
$H_1(M)$ is the  standard Sobolev space consisting of functions in
$L^2$ with gradient in $L^2$), an extremal function for the sharp
$L^2-$Nash inequality (\ref{E1.3}), that is, a function such that:
\begin{equation}
 \left( {\int_M {u_0^2 dV_g } }
\right)^{1 + \frac{2} {n}}\! \!= \!\left( {A_0(n)\!\int_M
\!{\left| {\nabla u_0} \right|_g^2 } dV_g  + B_{opt}(M)\!\int_M
{u_0^2 dV_g } } \right)\!\!\left( {\int_M {\left| u_0 \right|dV_g
} } \right)^{\frac{4} {n}}\nonumber
\end{equation}

In this paper we are focusing our interest in the special case
where the manifold is the standard unit sphere $\mathbb{S}^n$ of
$\mathbb{R}^{n+1}$. We study both Nash's inequalities  $L^1$ and
$L^2$ first in the general case and second in the presence of
symmetries.
\smallbreak
More precisely:
 \smallbreak
\noindent$\bullet\;$ We give the proof of the problem of finding
the first constant in the $L^2-$Nash inequality in $\mathbb{S}^n$
and we compute the exact value of the second best
constant of this inequality.\\
$\bullet\;$ We answer the problem of finding both best
constants in the $L^1-$Nash inequality in $\mathbb{S}^n$.\\
$\bullet\;$ We prove the existence of extremal functions in
$L^2$ and non existence in $L^1-$Nash inequalities.\\
$\bullet\;$ We study the problem of the best constants in the
$L^2-$Nash inequality in $\mathbb{S}^n$, $n\geq3$, where the data
are $G-$invariant under the action of the group $G=O(k)\times
O(m)$, $k+m=n+1$, $k\geq m\geq2$ and we find the best constants in
this case.

\section{Statement of results }
\begin{theorem}\label{T2.1}
For all  $\phi \in H_1(\mathbb{S}^n),\;n\geq 1$, there exists a
constant $B$ such that the following inequality holds
\begin{eqnarray}\label{E2.1}
\left(\int_{\mathbb{S}^n} {\phi^2 ds}\right)^{1+\frac{2}{n}}
\leqslant\left( A_0(n) {\int_{\mathbb{S}^n} {\left| {\nabla \phi}
\right|^2 ds } + B\int_{\mathbb{S}^n} {\phi^2 ds} } \right)
\left({\int_{\mathbb{S}^n} {\left| \phi \right|} ds }
\right)^{\frac{4} {{n }}}
\end{eqnarray}
Moreover the constant $A_0(n)$ is the optimal for this inequality.\\
\end{theorem}

\begin{theorem}\label{T2.2}
For all  $\phi \in H_1(\mathbb{S}^n),\;n\geq 1$, there exists a
constant $A$ such that the following inequality holds
\begin{eqnarray}\label{E2.2}
\left(\int_{\mathbb{S}^n} {\phi^2 ds}\right)^{1+\frac{2}{n}}
\leqslant \left( A{\int_{\mathbb{S}^n} {\left| {\nabla \phi}
\right|^2 ds } + \omega_n^{-\frac{2}{n}}\int_{\mathbb{S}^n}
{\phi^2 ds} } \right) \left({\int_{\mathbb{S}^n} {\left| \phi
\right|} ds } \right)^{\frac{4} {{n }}},
\end{eqnarray}
where $\omega_n$ denotes the volume of the standard unit sphere
$\mathbb{S}^n$ of $\mathbb{R}^{n+1}$. In particular
$$
\omega_{2n}=\frac{(4\pi)^n(n-1)!}{(2n-1)!}\;\;
 \mathrm{and}\;\;\omega_{2n+1}=\frac{2\pi^{n+1}}{n!}
 $$
Moreover $\omega_n^{-\frac{2}{n}}$ is the optimal constant for this inequality.\\
In addition there exists $\phi_0 \in
H_1(\mathbb{S}^n),\;\phi_0\not \equiv 0$, an  extremal function
for the sharp $L^2-$inequality
$(N(A_0(n),\omega_n^{-\frac{2}{n}})$, that is, such that
\begin{eqnarray}\label{E2.3}
\left(\int_{\mathbb{S}^n} {\phi^2_0 ds}\right)^{1+\frac{2}{n}} =
\left(A_0(n) {\int_{\mathbb{S}^n} {\left| {\nabla \phi_0}
\right|^2 ds } + \omega_n^{-\frac{2}{n}}\int_{\mathbb{S}^n}
{\phi^2_0 ds} } \right) \left({\int_{\mathbb{S}^n} {\left| \phi_0
\right|} ds } \right)^{\frac{4} {{n }}}
\end{eqnarray}
\end{theorem}

\begin{theorem}\label{T2.3}
For all  $\phi \in H_1(\mathbb{S}^n)$ there exists a constant
$B_\varepsilon$ such that the following inequality holds
\begin{eqnarray}\label{E2.4}
\left(\int_{\mathbb{S}^n} {\phi^2 ds}\right)^{1+\frac{2}{n}}
\leqslant \left(A_0(n)+\varepsilon\right){\int_{\mathbb{S}^n}
{\left| {\nabla \phi} \right|^2 ds } \left({\int_{\mathbb{S}^n}
{\left| \phi \right|} ds } \right)^{\frac{4} {{n }}}+
B_\varepsilon\left({\int_{\mathbb{S}^n} {\left| \phi \right|} ds }
\right)^{2+\frac{4} {{n }}}}
\end{eqnarray}
Moreover the constant $A_0(n)$
is the optimal for this inequality.\\
\end{theorem}

\begin{theorem}\label{T2.4}
For all  $\phi \in H_1(\mathbb{S}^n)$ there exists a constant $A$
such that the following inequality holds
\begin{eqnarray}\label{E2.5}
\left(\int_{\mathbb{S}^n} {\phi^2 ds}\right)^{1+\frac{2}{n}}
\leqslant A{\int_{\mathbb{S}^n} {\left| {\nabla \phi} \right|^2 ds
} \left({\int_{\mathbb{S}^n} {\left| \phi \right|} ds }
\right)^{\frac{4} {{n }}}+
\omega_n^{-1-\frac{2}{n}}\left({\int_{\mathbb{S}^n} {\left| \phi
\right|} ds } \right)^{2+\frac{4} {{n }}}}
\end{eqnarray}
Moreover $\omega_n^{-1-\frac{2}{n}}$ is the optimal constant for
this inequality.
\end{theorem}
\begin{corollary}\label{C2.1}
The inequality of Theorem \ref{T2.3} is false if $\varepsilon=0$.
\end{corollary}
\begin{corollary}\label{C2.2}
There do not exist extremal functions for the sharp $L^1-$Nash
inequality $N(A_0(n), B_{opt}^1)$.
\end{corollary}

\begin{theorem}\label{T2.5}
For all  $f \in H_{1,G}(\mathbb{S}^n),\;n\geq 3$, there exists a
constant $B$ such that the following inequality holds
\begin{eqnarray}\label{E2.1}
\left(\int_{\mathbb{S}^n} {f^2 ds}\right)^{1+\frac{2}{k}}
\leqslant\left(A_0(k)\omega_{n-k}^{-\frac{2}{k}}{\int_{\mathbb{S}^n}
{\left| {\nabla f} \right|^2 ds } + B\int_{\mathbb{S}^n} {f^2 ds}
} \right) \left({\int_{\mathbb{S}^n} {\left| f \right|} ds }
\right)^{\frac{4} {{k }}}
\end{eqnarray}
Moreover the constant $A_0(k)\omega_{n-k}^{-\frac{2}{k}}$ is the
optimal for this inequality.
\end{theorem}

\begin{theorem}\label{T2.6}
For all  $f\in H_{1,G}(\mathbb{S}^n),\;n\geq 3$, there exists a
constant $A$ such that the following inequality holds
\begin{eqnarray}\label{E2.2}
\left(\int_{\mathbb{S}^n} {f^2 ds}\right)^{1+\frac{2}{k}}
\leqslant \left( A{\int_{\mathbb{S}^n} {\left| {\nabla f}
\right|^2 ds } + \omega_n^{-\frac{2}{k}}\int_{\mathbb{S}^n} {f^2
ds} } \right) \left({\int_{\mathbb{S}^n} {\left| f \right|} ds }
\right)^{\frac{4} {{k }}}
\end{eqnarray}
Moreover the constant $\omega_n^{-\frac{2}{k}}$ is the optimal for
this inequality.
\end{theorem}

\section{Notations and preliminary results}
\subsection{The General Case}
Consider the sphere $\mathbb{S}^n \subset \mathbb{R}^{n+1}$, of
dimension $n$ and radius $1$. That is
$$
\mathbb{S}^n= \{x\in \mathbb{R}^{n+1}: |x| = 1\}
$$
The stereographic projection
$$
\Pi: \mathbb{S}^n\backslash\{N\} \rightarrow \mathbb{R}^n=\{x\in
\mathbb{R}^{n+1} : x_{n+1} = 0\}
$$
maps a point $P'\in \mathbb{S}^n\backslash\{N\}$ into the
intersection $P\in \mathbb{R}^n$ of the line joining $P'$ and the
north pole $N = (0,0, . . . , 1)$ with $\mathbb{R}^n$. \\
Let $g_{\alpha\beta}$ the standard metric of $\mathbb{S}^n$ (i.e.
the one inherited from $\mathbb{R}^{n+1}$) is expressed in terms
of stereographic coordinates by
$$
g_{\alpha\beta}=\left({\frac{2}{1+|x|^2}}\right)^2\delta_{\alpha\beta}.
$$
Hence the standard volume element of $\mathbb{S}^n$ is
$$
ds=\left({\frac{2}{1+|x|^2}}\right)^n dx$$

Let $H_{1}(\mathbb{S}^n)$ be the standard Sobolev space consisting
of functions in $L^2(\mathbb{S}^n)$ with gradient in
$L^2(\mathbb{S}^n)$. For any function $\phi \in
H_{1}(\mathbb{S}^n)$ set $u=\phi\circ\Pi^{-1}$. The integral and
the gradient Dirichlet integral, corresponding to a conformal
metric $ds = p^n dx$, where $p= {\frac{2}{1+|x|^2}} $, are:
\begin{eqnarray}\label{E3.1}
\int_{\mathbb{S}^n} \phi\, ds &=&\int_{\mathbb{R}^n} u\, p^n
dx\;\;\;\;
\end{eqnarray}
\begin{eqnarray}\label{E3.2}
\int_{\mathbb{S}^n}|\nabla \phi\,|^2
ds&=&\int_{\mathbb{R}^n}|\nabla u\,|^2 p^{n-2}dx
\end{eqnarray}

We may assume that $\mathbb{S}^n$ is covered by a finite number of
charts, say $(U_j, \xi_j), 1 \leq j\leq N$, such for any
$\varepsilon>0$, $\left( U_j, \xi_j\right)$ can be chosen such
that:
\begin{eqnarray}\label{E3.3}
1 - \varepsilon \leq \sqrt {\det ( {g_{\alpha\beta}^{j} } )}  \leq
1 + \varepsilon \quad \mathrm{on}\quad  U_j ,\quad
\mathrm{for}\quad 1 \leq \alpha,\beta \leq n
\end{eqnarray}
where the $g_{\alpha\beta}^j$'s are the components of $g$ in
$(U_j, \xi_j)$.\\
For each $j$ we consider $h_j \in C_{0}^\infty \left(\mathbb{R}^n
\right)$, $h_j\geq 0$ and set
\begin{eqnarray}\label{E3.4}
\eta _j  = \frac{{h_j \circ \xi _j }} {{\sum\nolimits_{j = 1}^N
{\left( {h_j \circ \xi _j } \right)} }}
\end{eqnarray}
The $\eta_j$'s are then a partition of unity for $\mathbb{S}^n$
relative to $U_j$'s.\\

\begin{lemma}\label{L3.1}
For any $\varepsilon >0$ and for all $\phi \in
C^\infty_0(\mathbb{S}^n)$ the following inequality holds
\begin{eqnarray}\label{E3.5}
\left( \int_{\mathbb{S}^n} {\left( {\eta _j \phi} \right)^2 ds}
\right)^{1 + \frac{2} {n}}\leqslant \left( A_0(n)+\varepsilon
\right){\int_{\mathbb{S}^n}{\left| {\nabla( \eta_j\phi)} \right|^2
} ds } \left( {\int_{\mathbb{S}^n} {\left| \eta_j\phi \right|ds }
} \right)^{\frac{4} {n}}
\end{eqnarray}
\end{lemma}
\textbf{Proof.}$\,\,\,$ By  (\ref{E3.1}) and  (\ref{E3.2}) because
of (\ref{E3.3}) for any  $ \phi \in C_0^\infty  ( \mathbb{S}^n ) $
and any $q\geq 1$ real, setting $(\eta _j \phi)\circ \Pi^{-1}=u_j
$ we obtain
\begin{equation}\label{E3.6}
\left( {1 - \varepsilon} \right)^n \int_{\mathbb{R}^n} {\left( {u
_j} \right)^q dx }\leq \int_{\mathbb{S}^n} {\left( {\eta _j \phi}
\right)^q ds} \leqslant \left( {1 + \varepsilon} \right)^n
\int_{\mathbb{R}^n} {\left( {u _j} \right)^q dx }
\end{equation}
and
\begin{equation}\label{E3.7}
\left( {1 - \varepsilon} \right)^{n-2} \int_{\mathbb{R}^n} {\left|
\nabla {u _j} \right|^2 dx }\leq \int_{\mathbb{S}^n} {\left|\nabla
(\eta _j \phi) \right|^2 ds} \leqslant \left( {1 + \varepsilon}
\right)^{n-2} \int_{\mathbb{R}^n} {\left| \nabla{u _j} \right|^2
dx }
\end{equation}
It is known, by Carlen and Loss \cite{Car-Los}, that for any $u\in
C_0^\infty (\mathbb{R}^n)$, the following inequality holds
\begin{eqnarray*}
 \left( {\int_{\mathbb{R}^n} {u^2 dx } }
\right)^{1 + \frac{2} {n}} \leqslant A_0(n)
{\int_{\mathbb{R}^n}{\left| {\nabla u} \right|^2 } dx   } \left(
{\int_{\mathbb{R}^n} {\left| u \right|dx } } \right)^{\frac{4}
{n}}
\end{eqnarray*}
For any $\varepsilon >0$,  we can choose $\delta>0$ such that for
any $x = \xi_j( s ) \in \mathbb{R}^n, s \in
U_j\subset\mathbb{S}^n$ and for all $u \in C_0^\infty  ( B_x (
\delta ) )$, ($B_x( \delta )\subset \xi_j( U_j )$ is the
$n-$dimensional ball of radius $\delta$ centered on $x$), the
following inequality holds
\begin{equation}\label{E3.8}
 \left( {\int_{\mathbb{R}^n} {u^2 dx } }
\right)^{1 + \frac{2} {n}} \leqslant  \left(
A_0(n)+\frac{\varepsilon}{2} \right){\int_{\mathbb{R}^n}{\left|
{\nabla u} \right|^2 } dx   } \left( {\int_{\mathbb{R}^n} {\left|
u \right|dx } } \right)^{\frac{4} {n}}
\end{equation}
From (\ref{E3.7})  because of (\ref{E3.5}) and (\ref{E3.6}) we
obtain
\begin{eqnarray}\label{E3.9}
\left( \int_{\mathbb{S}^n} {\left( {\eta _j \phi} \right)^2 ds}
\right)^{1 + \frac{2} {n}}&\leqslant & \left( {1 + \varepsilon}
\right)^{n+2} \left(\int_{\mathbb{R}^n} {\left( {u _j} \right)^2
dx }\right)^{1 + \frac{2} {n}}\nonumber\\
&\leqslant & \left( {1 + \varepsilon} \right)^{n+2}
A_0(n){\int_{\mathbb{R}^n}{\left| {\nabla u_j} \right|^2 } dx   }
\left( {\int_{\mathbb{R}^n} {\left| u_j \right|dx } }
\right)^{\frac{4} {n}}\nonumber\\
&\leqslant & \left( {1 + \varepsilon} \right)^{n+2}\frac{1}{\left(
{1 - \varepsilon} \right)^{n-2}} \frac{1}{\left( {1 + \varepsilon}
\right)^{4}}\nonumber\\
&&\times \left( A_0(n)+\frac{\varepsilon}{2}
\right){\int_{\mathbb{S}^n}{\left| {\nabla (\eta_j\phi)} \right|^2
} ds } \left( {\int_{\mathbb{S}^n} {\left| \eta_j\phi \right|ds }
}
\right)^{\frac{4} {n}}\nonumber\\
&=& \left( \frac{1+\varepsilon}{1-\varepsilon} \right)^{n-2}\left(
A_0(n)+\frac{\varepsilon}{2} \right){\int_{\mathbb{S}^n}{\left|
{\nabla( \eta_j\phi)} \right|^2 } ds }\nonumber\\
&&\times \left( {\int_{\mathbb{S}^n} {\left| \eta_j\phi \right|ds
} } \right)^{\frac{4} {n}}
\end{eqnarray}
Since the function $f:(0,1) \to (1, + \infty )$ with $f\left(
\varepsilon  \right) = \left( \frac{1+\varepsilon}{1-\varepsilon}
\right)^{n-2}$ is monotonically increasing, we can choose the
$\varepsilon>0$ such that the inequality
$$
f\left( \varepsilon  \right)\left( {A_0 ( n ) + \frac{\varepsilon
}{2}} \right) \le  A_0 ( n) + \varepsilon
$$
holds. Hence from (\ref{E3.9}) follows (\ref{E3.5}) and the lemma
is proved.\mbox{ }\hfill $\Box$

\begin{lemma}\label{L3.2}
For any $\varepsilon >0$ and for all $\phi \in
C^\infty_0(\mathbb{S}^n)$ there exists a constant
$B_\varepsilon>0$ such that the following inequality holds
\begin{eqnarray}\label{E3.10}
\left(\int_{\mathbb{S}^n} {\phi^2 ds}\right)^{1+\frac{2}{n}}
\leqslant \left(A_0(n)+\varepsilon\right)\left(
{\int_{\mathbb{S}^n} {\left| {\nabla \phi} \right|^2 ds } +
B_\varepsilon\int_{\mathbb{S}^n} {\phi^2 ds} } \right)
\left({\int_{\mathbb{S}^n} {\left| \phi \right|} ds }
\right)^{\frac{4} {{n }}}
\end{eqnarray}
\end{lemma}
\textbf{Proof.}$\,\,\,$ We set $\alpha _j = \frac{{\eta _j^2 }}
{{\sum\nolimits_{m = 1}^N {\eta _m^2 } }}, j = 1,2,..,N$, where
$\eta_j$ is defined by (\ref{E3.4}), and so $ \{ {\alpha _j }
\}_{j = 1,2,...N} $ is a partition of unity for $\mathbb{S}^n$
subordinated in the covering $( {U_j })_{j = 1,2,...,N}$,
functions $\sqrt {\alpha _j } $ are smooth and there exist a
positive  constant $H$ such that for any $j = 1,...,N$ holds
\begin{equation}\label{E3.11}
| {\nabla \sqrt {\alpha _j } } | \leqslant H
\end{equation}
Let  $\phi \in C^\infty(\mathbb{S}^n)$. Then we have
\begin{eqnarray}\label{E3.12}
\|\phi \|^2_{2}=\|\phi^2 \|_{1}= \|\sum\limits_{j =
1}^N\alpha_j\phi^2 \|_{1} \leq \sum\limits_{j =
1}^N\left\|\alpha_j\phi^2 \right\|_{1}=\sum\limits_{j =
1}^N\|\sqrt{\alpha_j}\phi \|^2_{2}
\end{eqnarray}
By Lemma \ref{L3.1}, for any $j$,
\begin{equation}\label{E3.13}
\|\sqrt{\alpha_j}\phi \|^2_{2}\leq
A^{\frac{n}{n+2}}\left\|\nabla(\sqrt{\alpha_j}\phi)\right\|^{\frac{2n}{n+2}}_{2}
\left\|\sqrt{\alpha_j}\phi\right\|^{\frac{4}{n+2}}_{1}
\end{equation}
where
$$
A=A_0(n)+\varepsilon.
$$
By H\"older's inequality,
\begin{equation}\label{E3.14}
\|\sqrt{\alpha_j}\phi \|_{1}\leq \|\alpha_j\phi
\|^{\frac{1}{2}}_{1}\,\, \|\phi \|^{\frac{1}{2}}_{1}
\end{equation}
As a consequence, by (\ref{E3.12}), (\ref{E3.13}) and
(\ref{E3.14}), for any $\phi \in C^\infty(\mathbb{S}^n)$ we obtain
\begin{eqnarray}\label{E4.11}
\int_{\mathbb{S}^n} {\phi^2 ds  }&\leqslant & A^{\frac{n}{n+2}}
\left( {\int_{\mathbb{S}^n} {\left| \phi \right|} ds }
\right)^{\frac{2} {{n+2}}}\sum\limits_{j = 1}^N {\left(
{\int_{\mathbb{S}^n} {\left| {\nabla \left( {\sqrt {\alpha _j }
\phi} \right)} \right|^2 ds } } \right)^{\frac{{n}} {{n+2}}} }
\nonumber\\
&&\times\left( {\int_{\mathbb{S}^n} {\alpha _j \left| \phi
\right|} ds } \right)^{\frac{2} {{n +2}}}
\end{eqnarray}
Moreover, for any $a_j ,b_j $ non negative and for all $p
\geqslant 1,q \geqslant 1$, with $ {\frac{1}{p}}  +  {\frac{1}{q}}
= 1$, by the H\"older's inequality in the discreet case it holds
\begin{equation}\label{E4.12}
\sum\limits_{j = 1}^N {a_j b_j } \leqslant \left( {\sum\limits_{j
= 1}^N {a_j^p } } \right)^{\frac{1}{p}} \left( {\sum\limits_{j =
1}^N {b_j^q } } \right)^{\frac{1}{q}}
\end{equation}
Setting  in (\ref{E4.12})
$$
a_j  = \left( {\int_{\mathbb{S}^n} {| {\nabla ( {\sqrt {\alpha _j
}\phi} )} |^2 ds } }\right)^{\frac{{n }} {{n +2}}}, \quad b_j  =
\left( {\int_{\mathbb{S}^n} {\alpha _j | \phi |} ds}
\right)^{\frac{2} {{n+2}}}
$$
$$
p = \frac{n+2} {n},\quad q =\frac{n+2}{2}
$$
we obtain
\begin{eqnarray}\label{E4.13}
&& \sum\limits_{j = 1}^N {\left( {\int_{\mathbb{S}^n} {\left|
{\nabla \left( {\sqrt {\alpha _j } \phi} \right)} \right|^2 ds } }
\right)^{\frac{{n}} {{n+2}}} } \left( {\int_{\mathbb{S}^n} {\alpha
_j \left| \phi \right|} ds } \right)^{\frac{2} {{n +2}}}\nonumber \hfill \\
&\leqslant& \left( {\sum\limits_{j = 1}^N {\int_{\mathbb{S}^n}
{\left| {\nabla \left( {\sqrt {\alpha _j }\phi} \right)} \right|^2
ds } } } \right)^{\frac{n} {{n +2}}} \left( {\sum\limits_{j = 1}^N
{\int_{\mathbb{S}^n}  {\alpha _j \left| \phi
\right|} ds } } \right)^{\frac{2} {{n +2}}} \nonumber\hfill \\
& = &\left( {\sum\limits_{j = 1}^N {\int_{\mathbb{S}^n}  {\left|
{\nabla \left( {\sqrt {\alpha _j } \phi} \right)} \right|^2 ds} }
} \right)^{\frac{n} {{n +2}}} \left( {\int_{\mathbb{S}^n} {\left(
{\sum\limits_{j = 1}^N {\alpha _j } } \right)\left| \phi \right|ds
} } \right)^{\frac{2} {{n +2}}}\nonumber \hfill \\ & =&
   \left( {\sum\limits_{j = 1}^N {\int_{\mathbb{S}^n} {\left|
{\nabla \left( {\sqrt {\alpha _j } \phi} \right)} \right|^2 ds } }
} \right)^{\frac{n} {{n +2}}}\left( {\int_{\mathbb{S}^n}
{\left|\phi \right|ds } } \right)^{\frac{2} {n+2}}
\end{eqnarray}
By (\ref{E4.11}) and (\ref{E4.13}) we obtain
\begin{eqnarray}\label{E4.14}
\int_{\mathbb{S}^n} {\phi^2 ds } &\leqslant & A^{\frac{n} {n+2}}
 \left( {\sum\limits_{j = 1}^N
{\int_{\mathbb{S}^n} {\left| {\nabla \left( {\sqrt {\alpha _j }
\phi} \right)} \right|^2 ds} } } \right)^{\frac{n } {n +2}}\left(
{\int_{\mathbb{S}^n} {\left| \phi\right|} ds } \right)^{\frac{4}
{n+2}}
\end{eqnarray}
Furthermore, since
\[
\left| {\nabla \left( {\sqrt {\alpha _j } \phi} \right)} \right|^2
= \alpha _j \left| {\nabla \phi} \right|^2  + \phi^2 \left|
{\nabla \left( {\sqrt {\alpha _j } } \right)} \right|^2 +
2\left\langle {\nabla\phi,\nabla \left( {\sqrt {\alpha _j } }
\right)} \right\rangle \phi\sqrt {\alpha _j }
\]
and because of  (\ref{E3.11}), we obtain
\begin{eqnarray*}
  \sum\limits_{j = 1}^N {\int_{\mathbb{S}^n} {\left| {\nabla
   \left( {\sqrt {\alpha _j } \phi} \right)} \right|^2 ds}
    }  &=& \sum\limits_{j = 1}^N {\int_{\mathbb{S}^n} {\alpha _j
    \left| {\nabla \phi} \right|^2 ds + \sum\limits_{j =
    1}^N {\int_{\mathbb{S}^n}{\left| {\nabla \left( {\sqrt
    {\alpha _j } } \right)} \right|^2 \phi^2 ds } } } }  \hfill \\
  && + 2\sum\limits_{j = 1}^N {\int_{\mathbb{S}^n} {\left\langle
   {\nabla \left( {\sqrt {\alpha _j } } \right),\nabla \phi}
   \right\rangle \phi\sqrt {\alpha _j } ds } }  \hfill \\
 \end{eqnarray*}
\begin{eqnarray*} &  =& \sum\limits_{j = 1}^N {\int_{\mathbb{S}^n} {\alpha _j
  \left| {\nabla \phi} \right|^2 ds + \sum\limits_{j = 1}^N
   {\int_{\mathbb{S}^n} {\left| {\nabla \left( {\sqrt {\alpha _j }
    } \right)} \right|^2 \phi^2 ds} } } }  \hfill \\
  && + 2\int_{\mathbb{S}^n} {\sum\limits_{j = 1}^N
  {\left\langle {\nabla \left( {\sqrt {\alpha _j }
  } \right),\nabla \phi} \right\rangle\phi\sqrt
  {\alpha _j } ds } }  \hfill \\
  &=& \int_{\mathbb{S}^n} {\left({\sum\limits_{j = 1}^N {\alpha _j } } \right)\left|
  {\nabla \phi} \right|^2 ds} \\&& + \int_{\mathbb{S}^n}
  {\sum\limits_{j = 1}^N {\alpha _j } \left|
  {\nabla \left( {\sqrt {\alpha _j } } \right)}
  \right|^2 \phi^2 ds }  \hfill \\
  &&+ 2\int_{\mathbb{S}^n}
  {\sum\limits_{j = 1}^N {\left\langle {\nabla
  \left( {\sqrt {\alpha _j } } \right),\nabla \phi}
  \right\rangle \phi\sqrt {\alpha _j } ds } }  \hfill \\
 & \leqslant & \int_{\mathbb{S}^n} {\left|
  {\nabla \phi} \right|^2 ds }  + HN\int_{\mathbb{S}^n}
  {\phi^2 ds }  \hfill \\
 && + 2\int_{\mathbb{S}^n}
  {\sum\limits_{j = 1}^N {\left\langle {\nabla
   \left( {\sqrt {\alpha _j } } \right),\nabla \phi}
   \right\rangle \phi\sqrt {\alpha _j } ds } }  \hfill \\
 & =& \int_{\mathbb{S}^n} {\left|
  {\nabla \phi} \right|^2 ds }  + C_\varepsilon\int_{\mathbb{S}^n} {\phi^2 ds }  \hfill \\
  &&+ 2\int_{\mathbb{S}^n}
   {\sum\limits_{j = 1}^N {\left\langle {\nabla
    \left( {\sqrt {\alpha _j } } \right),\nabla \phi}
    \right\rangle \phi\sqrt {\alpha _j } ds } }  \hfill \\
\end{eqnarray*}
But
\[
2\sum\limits_{j = 1}^N {\sqrt {\alpha _j } \nabla \left( {\sqrt
{\alpha _j } } \right) = } \sum\limits_{j = 1}^N {2\sqrt {\alpha
_j } \nabla \left( {\sqrt {\alpha _j } } \right) = }
\sum\limits_{j = 1}^N {\left( {\nabla \alpha _j } \right) = }
\nabla ( {\sum\limits_{j = 1}^N {\alpha _j } } ) = 0
\]
Thus the following inequality is true
\begin{equation}\label{E4.15}
\sum\limits_{j = 1}^N {\int_{\mathbb{S}^n}  {\left| {\nabla \left(
{\sqrt {\alpha _j } \phi} \right)} \right|^2 ds } }  \leqslant
\int_{\mathbb{S}^n}  {\left| {\nabla \phi} \right|^2 ds }  +
B_\varepsilon\int_{\mathbb{S}^n}  {\phi^2 ds }
\end{equation}
Hence by (\ref{E4.14}) and (\ref{E4.15}) we obtain
\[
\begin{gathered}
\int_{\mathbb{S}^n} {\phi^2 ds} \leqslant A^{\frac{n} {n+2}}
 \left( {\int_{\mathbb{S}^n} {\left|
{\nabla \phi} \right|^2 ds } + B_\varepsilon\int_{\mathbb{S}^n}
{\phi^2 ds} } \right)^{\frac{{n}} {{n+2}}}
\left({\int_{\mathbb{S}^n} {\left| \phi \right|} ds }
\right)^{\frac{4} {{n +2}}} \hfill \\
\end{gathered}
\]
or
\[
\begin{gathered}
\left(\int_{\mathbb{S}^n} {\phi^2 ds}\right)^{1+\frac{2}{n}}
\leqslant \left(A_0(n)+\varepsilon\right)\left(
{\int_{\mathbb{S}^n} {\left| {\nabla \phi} \right|^2 ds } +
B_\varepsilon\int_{\mathbb{S}^n} {\phi^2 ds} } \right)
\left({\int_{\mathbb{S}^n} {\left| \phi \right|} ds }
\right)^{\frac{4} {{n }}} \hfill \\
\end{gathered}
\]
and the lemma is proved.\mbox{ }\hfill $\Box$
\subsection{The Case of Existence of Symmetries for $n\geq3$}
Let
$$
\mathbb{R}^{n+1}=\mathbb{R}^k\times\mathbb{R}^m=\{(x,y):x\in\mathbb{R}^k,\;\;
y\in\mathbb{R}^m\},
$$
where $k+m=n+1,\quad k\geq m\geq 2 $.\\
Then
$$
\mathbb{S}^n=\{(x,y):|x|^2+|y|^2=1\}
$$
Let $x=(x^1, x^2,...,x^k)\in \mathbb{R}^k$ and $y=(x^{k+1},
x^{k+2},...,x^{n+1})\in \mathbb{R}^m$, where $\{x^i,\,
i=1,2,...,n+1\}$ is a coordinate system of $\mathbb{R}^{n+1}$.

It is well known that $\mathbb{S}^n$ enjoys a lot of symmetries,
namely, the compact Lie group $O(n +1)$ acts isometrically on
$\mathbb{S}^n$. Let now $G=O(k)\times O(m)$. Then $G$ is a compact
subgroup of $O(n +1)$. For $g=(g_1,g_2)\in G$, where $g_1\in O(k)$
and $g_2\in O(m)$, the action of $G$ on $S^n$ is defined by
$g(x,y)=(g_1 x,g_2 y)$ and if $P(x,y)\in\mathbb{S}^n$ its orbit
under the action of $G$, since $|x|^2+|y|^2=1$, is
$$
O_P=\mathbb{S}^{k-1}(|x|)\times
\mathbb{S}^{m-1}(|y|)=\mathbb{S}^{k-1}(|x|)\times
\mathbb{S}^{m-1}(\sqrt{1-|x|^2})
$$
Denote $C_G^\infty (\mathbb{\mathbb{S}}^n)$ the space of all
$G-$invariant functions under the action  of the group $G$ and
$H_{1, G}(\mathbb{S}^n)$ the space of all $G-$invariant functions
of $H_1(\mathbb{S}^n)$. Under the above considerations, if $f \in
H_{1, G}(\mathbb{S}^n)$, we can set $|x|= sin \theta\;$ and
$\;|y|= cos \theta$, $0\leq \theta \leq \pi/2$ and then $f$ is a
function of one variable $\theta$ and the following formulas hold:
\begin{equation}\label{E1}
\int_{\mathbb{S}^n}|f|^q ds=\omega_{k-1}\omega_{m-1}\int_0^{\pi/2}
|f|^q sin^{k-1}\theta\, \,cos^{m-1}\theta d\theta
\end{equation}
\begin{equation}\label{E2}
\int_{\mathbb{S}^n}|\nabla f|^2
ds=\omega_{k-1}\omega_{m-1}\int_0^{\pi/2} (f')^2 sin^{k-1}\theta\,
\,cos^{m-1}\theta d\theta
\end{equation}
\section{Proofs}
{\bf{Proof of Theorem \ref{T2.1}}.} For $n=1$, the theorem is
true, see Theorem 2.5 in \cite{Hum1}. Let $n\geq 2$. In order to
prove inequality (\ref{E2.1}) it is equivalence to proving  that
for all $\phi \in H^2_1(\mathbb{S}^n)$ there exists a constant
$B'$ such that the following inequality holds
\begin{eqnarray*}
\left(\int_{\mathbb{S}^n} {\phi^2 ds}\right)^{1+\frac{2}{n}}
\leqslant A_0(n)\left( {\int_{\mathbb{S}^n} {\left| {\nabla \phi}
\right|^2 ds } + B'\int_{\mathbb{S}^n} {\phi^2 ds} } \right)
\left({\int_{\mathbb{S}^n} {\left| \phi \right|} ds }
\right)^{\frac{4} {{n }}}
\end{eqnarray*}
We use a proof based on Lemma \ref{L3.2}. Suppose by contradiction
that the inequality is not true. Then for any $\alpha>0$ there
exists $\phi_\alpha   \in C_0^\infty \left(\mathbb{S}^n \right)$
such that
\begin{equation}\label{E4.1}
\frac{\left( {\int_{\mathbb{S}^n} {\left| {\nabla \phi_\alpha}
\right|^2 ds } +\alpha\int_{\mathbb{S}^n} {\phi_\alpha^2 ds} }
\right) \left({\int_{\mathbb{S}^n} {\left| \phi_\alpha \right|} ds
} \right)^{\frac{4} {{n }}}}{\left(\int_{\mathbb{S}^n}
{\phi_\alpha^2 ds}\right)^{1+\frac{2}{n}}} < \frac{1}{A_0(n)}
\end{equation}
By (\ref{E4.1}) because of (\ref{E3.1}) and (\ref{E3.2}) we obtain
equivalently
\begin{equation}\label{E4.2}
\frac{\left( {\int_{\mathbb{R}^n} {\left| {\nabla u_\alpha}
\right|^2 p^{n-2}dx } + \alpha\int_{\mathbb{R}^n} {u_\alpha^2 p^n
dx} } \right) \left({\int_{\mathbb{R}^n} {\left| u_\alpha
\right|}p^n dx } \right)^{\frac{4} {{n
}}}}{\left(\int_{\mathbb{R}^n} {u_\alpha^2p^n
dx}\right)^{1+\frac{2}{n}}} <  \frac{1}{A_0(n)}
\end{equation}
where $u_\alpha=\phi_\alpha\circ\Pi^{-1}$.

For any $\lambda >0$, define $u_{\alpha_{_\lambda}}$ by
$u_{\alpha_{_\lambda}}(x)=u_\alpha(\lambda x)$. So, for any
$\lambda$, $u_{\alpha_{_\lambda}}$ has compact support and since
$p=\frac{2}{1+|x|^2}$ the following hold
\begin{equation}\label{E4.3}
\int_{\mathbb{R}^n} \left| \nabla u_\alpha \right|^2
p^{n-2}dx=\lambda^{n-2} \int_{\mathbb{R}^n} | \nabla
u_{\alpha_{_\lambda}} |^2
\left(\frac{2}{1+|\frac{x}{\lambda}|^2}\right)^{n-2}dx
\end{equation}
\begin{equation}\label{E4.4}
\int_{\mathbb{R}^n}  u_\alpha^2 p^{n}dx=\lambda^{n}
\int_{\mathbb{R}^n}  u_{\alpha_{_\lambda}}^2
\left(\frac{2}{1+|\frac{x}{\lambda}|^2}\right)^{n}dx
\end{equation}
\begin{equation}\label{E4.5}
\left(\int_{\mathbb{R}^n} \left|  u_\alpha \right|
p^{n}dx\right)^\frac{4}{n}=\lambda^4 \left(\int_{\mathbb{R}^n} |
u_{\alpha_{_\lambda}} |
\left(\frac{2}{1+|\frac{x}{\lambda}|^2}\right)^{n}dx\right)^\frac{4}{n}
\end{equation}
\begin{equation}\label{E4.6}
\left(\int_{\mathbb{R}^n}  u_\alpha^2
p^{n}dx\right)^{1+\frac{2}{n}}=\lambda^{n+2}
\left(\int_{\mathbb{R}^n}  u_{\alpha_{_\lambda}}^2
\left(\frac{2}{1+|\frac{x}{\lambda}|^2}\right)^{n}dx\right)^{1+\frac{2}{n}}
\end{equation}
By (\ref{E4.2}) because of (\ref{E4.3}), (\ref{E4.4}),
(\ref{E4.5}) and (\ref{E4.6}), for $\lambda\rightarrow\infty $, we
obtain
\begin{eqnarray*}
\frac{  {\int_{\mathbb{R}^n} | \nabla u_{\alpha_{_\lambda}} |^2
2^{n-2}dx  }  \left(\int_{\mathbb{R}^n} | u_{\alpha_{_\lambda}} |
2^{n}dx\right)^\frac{4}{n}}{ \left(\int_{\mathbb{R}^n}
u_{\alpha_{_\lambda}}^2 2^{n}dx\right)^{1+\frac{2}{n}}} <
\frac{1}{A_0(n)}
\end{eqnarray*}
or
\begin{eqnarray*}
\frac{  {\int_{\mathbb{R}^n} | \nabla u_{\alpha_{_\lambda}} |^2
dx  }  \left(\int_{\mathbb{R}^n} | u_{\alpha_{_\lambda}} |
dx\right)^\frac{4}{n}}{ \left(\int_{\mathbb{R}^n}
u_{\alpha_{_\lambda}}^2 dx\right)^{1+\frac{2}{n}}} <
\frac{1}{A_0(n)}
\end{eqnarray*}
Because of Carlen-Loss  Theorem \cite{Car-Los}, last inequality is
false and the theorem is proved.\mbox{ }\hfill $\Box$\\

\noindent{\bf{Proof of Theorem \ref{T2.2}}.} For $n=1$,
the theorem is true, see Corollary 5.2 in \cite{Hum1}.\\
If $n=2$, by Theorem $1.2$ in \cite{Hum1}, we produce that
\begin{eqnarray}\label{E4.7}
B_{opt} {(\mathbb{S}^2)}\geq \omega_2^{-1}
\end{eqnarray}
Let $n\geq 3$. By Theorem \ref{T2.1} follows that for all  $\phi
\in H^2_1(\mathbb{S}^n)$ there exists a constant $B$ such that the
following inequality holds
\begin{eqnarray}\label{E4.8}
\left(\int_{\mathbb{S}^n} {\phi^2 ds}\right)^{1+\frac{2}{n}}
\leqslant \left(A_0(n) {\int_{\mathbb{S}^n} {\left| {\nabla \phi}
\right|^2 ds } + B\int_{\mathbb{S}^n} {\phi^2 ds} } \right)
\left({\int_{\mathbb{S}^n} {\left| \phi \right|} ds }
\right)^{\frac{4} {{n }}}
\end{eqnarray}
On the one hand, by taking $\phi=1$ in (\ref{E4.8}), one obtains
that $B\geq \omega_n^{-\frac{2}{n}}$. In particular
\begin{eqnarray}\label{E4.9}
B_{opt} {(\mathbb{S}^n)}\geq \omega_n^{-\frac{2}{n}}
\end{eqnarray}
On the other hand, by H\"older's inequality, for any $\phi \in
H^2_1(\mathbb{S}^n)$ and for $p=\frac{2n}{n-2}$, it holds that
\begin{eqnarray*}
\int_{\mathbb{S}^n} {\phi^2 ds} \leqslant
\left(\int_{\mathbb{S}^n} {\phi^p
ds}\right)^\frac{1}{p-1}\left(\int_{\mathbb{S}^n} {|\phi|
ds}\right)^\frac{p-2}{p-1}
\end{eqnarray*}
or this
\begin{eqnarray}\label{E4.10}
\left(\int_{\mathbb{S}^n} {\phi^2 ds}\right)^\frac{2(p-1)}{p}
\leqslant \left(\int_{\mathbb{S}^n} {\phi^p
ds}\right)^\frac{2}{p}\left(\int_{\mathbb{S}^n} {|\phi|
ds}\right)^\frac{2(p-2)}{p}
\end{eqnarray}
By Theorem $4.2$ in \cite{Heb}, there exists $A\in\mathbb{R}$ such
that for any $\phi \in H^2_1(\mathbb{S}^n)$, holds
\begin{eqnarray}\label{E4.11}
 \left(\int_{\mathbb{S}^n} {\phi^p
ds}\right)^\frac{2}{p}\leq A \int_{\mathbb{S}^n} {|\nabla
\phi|^2ds }+ \omega_n^{-\frac{2}{n}}\int_{\mathbb{S}^n}  \phi^2ds
\end{eqnarray}
By (\ref{E4.10}), because of (\ref{E4.11}), we have
\begin{eqnarray*}
\left(\int_{\mathbb{S}^n} {\phi^2 ds}\right)^\frac{2(p-1)}{p}
\leqslant \left(A \int_{\mathbb{S}^n} {|\nabla \phi|^2ds }+
\omega_n^{-\frac{2}{n}}\int_{\mathbb{S}^n}
\phi^2ds\right)\left(\int_{\mathbb{S}^n} {|\phi|
ds}\right)^\frac{2(p-2)}{p}
\end{eqnarray*}
and since
$$
\frac{2(p-1)}{p}=2-\frac{2}{p}=2-\frac{n-2}{n}=1+\frac{2}{n},
$$
$$
\frac{2(p-2)}{p}=2-\frac{4}{p}=2-\frac{2(n-2)}{n}=\frac{4}{n}
$$
we finally obtain
\begin{eqnarray*}
\left(\int_{\mathbb{S}^n} {\phi^2 ds}\right)^{1+\frac{2}{n}}
\leqslant \left(A \int_{\mathbb{S}^n} {|\nabla \phi|^2ds }+
\omega_n^{-\frac{2}{n}}\int_{\mathbb{S}^n}
\phi^2ds\right)\left(\int_{\mathbb{S}^n} {|\phi|
ds}\right)^\frac{4}{n}
\end{eqnarray*}
From this inequality and the definition of $B_{opt}
(\mathbb{S}^n)$, we obtain
\begin{eqnarray}\label{E4.12}
\omega_n^{-\frac{2}{n}}\geq B_{opt} {(\mathbb{S}^n)},\;\;
\forall\;\;n\geq2.
\end{eqnarray}
Further (\ref{E4.7}), (\ref{E4.9}) and (\ref{E4.12}) yield
\begin{eqnarray*}
B_{opt} {(\mathbb{S}^n)}= \omega_n^{-\frac{2}{n}},\;\;
\forall\;\;n\geq2
\end{eqnarray*}

For the second part of the Theorem, suppose by contradiction that
for all $\phi \in C_0^\infty \left(\mathbb{S}^n \right)$ the
following inequality holds
\begin{eqnarray*}
\left(\int_{\mathbb{S}^n} {\phi^2 ds}\right)^{1+\frac{2}{n}} <
A_0(n)\left( {\int_{\mathbb{S}^n} {\left| {\nabla \phi} \right|^2
ds } + A^{-1}_0(n)\omega_n^{-\frac{2}{n}}\int_{\mathbb{S}^n}
{\phi^2_0 ds} } \right) \left({\int_{\mathbb{S}^n} {\left| \phi
\right|} ds } \right)^{\frac{4} {{n }}}
\end{eqnarray*}
or
\begin{eqnarray*}
\frac{\left( {\int_{\mathbb{S}^n} {\left| {\nabla \phi} \right|^2
ds } +A^{-1}_0(n)\omega_n^{-\frac{2}{n}}\int_{\mathbb{S}^n}
{\phi^2 ds} } \right) \left({\int_{\mathbb{S}^n} {\left| \phi
\right|} ds } \right)^{\frac{4} {{n }}}}{\left(\int_{\mathbb{S}^n}
{\phi^2 ds}\right)^{1+\frac{2}{n}}} > \frac{1}{A_0(n)}
\end{eqnarray*}
Following the same steps as in the first part of theorem we
conclude that for all $u\in C_0^\infty \left(\mathbb{R}^n \right)$
\begin{eqnarray*}
\left(\int_{\mathbb{R}^n} u^2 dx\right)^{1+\frac{2}{n}} < A_0(n)
\int_{\mathbb{R}^n} | \nabla u |^2 dx  \left(\int_{\mathbb{R}^n} |
u| dx\right)^\frac{4}{n},
\end{eqnarray*}
which is false since, according to \cite{Car-Los}, there exists an
integrable  function $f_n$ on $\mathbb{R}^n$, such that its
distributional gradient is a square integrable function such that,
the equality bellow holds
\begin{eqnarray*}
\left(\int_{\mathbb{R}^n} f_n^2 dx\right)^{1+\frac{2}{n}} =A_0(n)
\int_{\mathbb{R}^n} | \nabla f_n |^2 dx  \left(\int_{\mathbb{R}^n}
| f_n| dx\right)^\frac{4}{n}.
\end{eqnarray*}
Moreover, it is easy to verify that constant functions are
extremal functions for the sharp $L^2-$Nash inequality and the
theorem is proved.\mbox{ }\hfill $\Box$\\

\noindent{\bf{Proof of Theorem \ref{T2.3}}.} Let
$$ a=\int_{\mathbb{S}^n} {\phi^2 ds} \quad \mathrm{and }\;\;
b=\left(\int_{\mathbb{S}^n} {|\phi| ds}\right)^\frac{4}{n} $$
Mimicking what is done in \cite{Dru-Heb-Vau}, let $\varepsilon_1
>0$ to be chosen later on, and set
$$ p=\frac{n+2}{n} \;\;\mathrm{and }\;\;q=\frac{n+2}{2} $$ Then $$
\frac{1}{p}+\frac{1}{q}=1$$ and so, by the elementary inequality
$$
xy\leq \frac{x^p}{p}+\frac{y^q}{q}\;\;\mathrm{for\;\; all }\;\;
 x,y\geq 0,\;\;\mathrm{and\;\;  for\;\;  all}\;\; \;\; p,q\geq 0\;\;
s.t.\;\; \frac{1}{p}+\frac{1}{q}=1,
$$
for $x=a\varepsilon_1$ and $y=\frac{b}{\varepsilon_1}$ we obtain
\begin{equation}\label{E4.13}
\int_{\mathbb{S}^n} {\phi^2 ds} \left(\int_{\mathbb{S}^n} {|\phi|
ds}\right)^\frac{4}{n}\leq
\frac{n{\varepsilon_1}^\frac{n+2}{n}}{n+2}\left(
\int_{\mathbb{S}^n} {\phi^2
ds}\right)^{1+\frac{2}{n}}+\frac{2\varepsilon_1^{-\frac{n+2}{2}}}{n+2}\left(\int_{\mathbb{S}^n}
{|\phi| ds}\right)^{2+\frac{4}{n}}
\end{equation}
By Lemma \ref{L3.2} arises that,  for any $\varepsilon >0$ and for
all $\phi \in C^\infty_0(\mathbb{S}^n)$ there exists a constant
$B_\varepsilon>0$ such that the following inequality holds
\begin{eqnarray}\label{E4.14}
\left(\int_{\mathbb{S}^n} {\phi^2 ds}\right)^{1+\frac{2}{n}}
\leqslant  A\left( {\int_{\mathbb{S}^n} {\left| {\nabla \phi}
\right|^2 ds } + B_\varepsilon\int_{\mathbb{S}^n} {\phi^2 ds} }
\right) \left({\int_{\mathbb{S}^n} {\left| \phi \right|} ds }
\right)^{\frac{4} {{n }}}
\end{eqnarray}
where
$$
A=A_0(n)+\frac{\varepsilon}{2}.
$$
Combining (\ref{E4.13}) and (\ref{E4.14}) we obtain
\begin{eqnarray*}
\left(\int_{\mathbb{S}^n} {\phi^2 ds}\right)^{1+\frac{2}{n}}
\leqslant \frac{A}{C}{\int_{\mathbb{S}^n} {\left| {\nabla \phi}
\right|^2 ds } \left({\int_{\mathbb{S}^n} {\left| \phi \right|} ds
} \right)^{\frac{4} {{n }}}+ B\left({\int_{\mathbb{S}^n} {\left|
\phi \right|} ds } \right)^{2+\frac{4} {{n }}}}
\end{eqnarray*}
where
$$
C=1- \frac{n}{n+2}{\varepsilon_1}^\frac{n+2}{n}AB_\varepsilon
\quad \mathrm{and }\quad
B=\frac{AB_\varepsilon}{C}\frac{2}{n+2}\varepsilon_1^{-\frac{n+2}{2}}
$$
We can choose $\varepsilon_1$ such that
$$\frac{A}{C}=A_0(n)+\varepsilon$$
and the theorem is proved.\\

\noindent Proof of Theorem \ref{T2.4} was discussed in
\cite{Dru-Heb-Vau}, (see Theorem 3.1).\mbox{ }\hfill $\Box$\\

\noindent{{\bf{Proof of Corollary  \ref{C2.1}}.} Since the scalar
curvature of $\mathbb{S}^n$ is $n(n-1)>0$ our result arises
immediately from Theorem $1.3$ of \cite{Dru-Heb-Vau}.\mbox{
}\hfill $\Box$
 \bigbreak \noindent{{\bf{Proof of Corollary
\ref{C2.1}}.} The conclusion arises immediately by Theorems
\ref{T2.3} and \ref{T2.4}.\mbox{ }\hfill $\Box$ \bigbreak
\noindent{\bf{Proof of Theorem \ref{T2.5}}.}  Let $\varepsilon>0$
be given. We consider $P \in M$ and its orbit $O_P $ of dimension
$k$. For any $Q = \tau(P) \in O_P $, where $\tau \in G$, we build
a chart around $Q$, denoted by $( {\tau ( \Omega_P )\!, \xi_P
\circ \tau ^{ - 1} } )$ and ``isometric'' to $( {\Omega_P ,\xi_P }
)$. $ O_P $ is then covered by such charts. We denote by $(
{\Omega _{m} } )_ {m = 1,...,M} $ a finite extract covering. We
then choose $\delta>0$ small enough, depending on $P$ and
$\varepsilon$, such that $ O_{\!P,\; \delta } = \left\{ {Q \in
\mathbb{S}^n: d(Q,O_P ) < \delta } \right\} $ the neighborhood $
O_{\!P,\; \delta }$, (where $d( { \cdot ,O_P }
)$ is the distance to the orbit) has the following properties:\\
$(i)\;$ $\overline { O_{\!P,\; \delta }} $ is a
submanifold of $\mathbb{S}^n$ with boundary,\\
$(ii)\;$ $d^2 ( { \cdot ,O_P } )$,  is a $C^\infty $
function on $ O_{\!P,\; \delta } $ and \\
$(iii)\;$ $ O_{\!P,\; \delta } $ is covered by $(
{\Omega _m } )_{m = 1,...,M} $.\\
Clearly, $\mathbb{S}^n$ is covered by $\cup_{P \in \mathbb{S}^n}
O_{\!P,\; \delta } $. We denote by $( {O_{j, \,\delta } } )_{j =
1,...,J} $ a finite extract covering of $\mathbb{S}^n$, where all
$O_{j, \,\delta } $'s are covered by $( {\Omega _{jm} } )_{m =
1,...,M_j } $. On each $( {O_{j, \,\delta } } ),{j = 1,...,J} $ we
consider functions depending only on the distance to $O_P$, and we
build a partition of unity $(\eta_j)$ relative to ${O_{j, \,\delta
} }$  such that for any $j$, $\eta_j\in C^\infty_G$. For any $f\in
C^\infty_G$, $\eta_j f\in C^\infty_G$ has compact support in
${O_{j, \,\delta } } $ and is a function of one variable. Thus
this partition of unity corresponds a subdivision of the interval
of integration $[0, \pi / 2]$ consisted of $J$ subintervals
$[\theta_{j-1}, \theta_j]$,  not
necessarily of equal length. \\
For any  subinterval $[\theta_{j-1}, \theta_j]$ there exists a
small $\varepsilon_j>0$, such that
\begin{eqnarray}\label{E3}
(1+\varepsilon_j)cos^{m-1}\theta_j =1
\end{eqnarray}

By Lemma \ref{L3.1} applied in $\mathbb{S}^k$ for any
$\varepsilon_0 >0$ and for all $\phi\in C_0^\infty (\mathbb{S}^k)$
the following inequality holds
\begin{eqnarray*}
\left(\int_{\mathbb{S}^k} {\phi^2 ds}\right)^{1+\frac{2}{k}}
\leqslant \left( A_0(k)+ \frac{\varepsilon_0}{2}\right)
\int_{\mathbb{S}^k} \left|\nabla \phi \right|^2 ds
\left(\int_{\mathbb{S}^k} \left| \phi \right| ds \right)^{\frac{4}
{{k }}},
\end{eqnarray*}
So for a radial function $\phi$ we obtain
\begin{eqnarray}\label{E4}
\left(\omega_{k-1}\int_0^{\pi/2}\!\! {\phi^2 sin^{k-1}\theta
d\theta}\right)^{1+\frac{2}{k}}\!\!\!& \leqslant&\!\!\left(
A_0(k)+\frac{\varepsilon_0}{2}
\right)\left(\omega_{k-1}\int_0^{\pi/2} \!\!\left( \phi' \right)^2
sin^{k-1}\theta d\theta\right)\nonumber\\&&\times
\left({\omega_{k-1}\int_0^{\pi/2}\!\!{\left|\phi \right|}
sin^{k-1}\theta d\theta } \right)^{\frac{4} {{k }}}
\end{eqnarray}

Let $f\in C^\infty_G$. Then $f$ is a function of one variable and
$\eta_j f\in C^\infty_G$ has compact support in ${O_{j, \,\delta }
} $ which corresponds to the  subinterval $[\theta_{j-1},
\theta_j]$.\\
By (\ref{E4}), because of (\ref{E3}),  we obtain
\begin{eqnarray*}
 \left(\omega_{k-1}\omega_{m-1}\int_{\theta_{j-1}}^{\theta_j} {(\eta_j f)^2
sin^{k-1}\theta cos^{m-1}\theta d\theta}\right)^{1+\frac{2}{k}}
\leq
\end{eqnarray*}
\begin{eqnarray*}
\left(\omega_{k-1}\omega_{m-1}\int_{\theta_{j-1}}^{\theta_j}
{(\eta_j f)^2 sin^{k-1}\theta d\theta}\right)^{1+\frac{2}{k}} \leq
\end{eqnarray*}
\begin{eqnarray*}
\left( A_0(k)+\frac{\varepsilon_0}{2}
\right)\omega_{m-1}^{-\frac{2}{k}}\left(\omega_{k-1}\omega_{m-1}\int_{\theta_{j-1}}^{\theta_j}
\!\!\left( (\eta_j f)' \right)^2 sin^{k-1}\theta
d\theta\right)\\\times
\left({\omega_{k-1}\omega_{m-1}\int_{\theta_{j-1}}^{\theta_j}\!\!{\left|\eta_j
f \right|} sin^{k-1}\theta d\theta } \right)^{\frac{4} {{k }}}
\leq
\end{eqnarray*}
\begin{eqnarray*}
\left( A_0(k)+\frac{\varepsilon_0 }{2}\right)(1+\varepsilon_j)
\omega_{m-1}^{-\frac{2}{k}}\left(\omega_{k-1}\omega_{m-1}\int_{\theta_{j-1}}^{\theta_j}
\!\!\left( (\eta_j f)' \right)^2 sin^{k-1}\theta cos^{m-1}\theta
d\theta\right)\\\times
\left({\omega_{k-1}\omega_{m-1}\int_{\theta_{j-1}}^{\theta_j}\!\!{\left|\eta_j
f \right|} sin^{k-1}\theta cos^{m-1}\theta d\theta }
\right)^{\frac{4} {{k }}}
\end{eqnarray*}
or
\begin{eqnarray*}
 \left(\omega_{k-1}\omega_{m-1}\int_0^{\pi/2} {(\eta_j f)^2
sin^{k-1}\theta cos^{m-1}\theta d\theta}\right)^{1+\frac{2}{k}}
\leq
\end{eqnarray*}
\begin{eqnarray}\label{E5}
\nonumber \\\left( A_0(k)+\frac{\varepsilon_0}{2}
\right)(1+\varepsilon_j)
\omega_{m-1}^{-\frac{2}{k}}\left(\omega_{k-1}\omega_{m-1}\int_0^{\pi/2}
\!\!\left( (\eta_j f)' \right)^2 sin^{k-1}\theta cos^{m-1}\theta
d\theta\right)\nonumber\\\times
\left({\omega_{k-1}\omega_{m-1}\int_0^{\pi/2}\!\!{\left|\eta_j f
\right|} sin^{k-1}\theta cos^{m-1}\theta d\theta }
\right)^{\frac{4} {{k }}}\nonumber\\
\end{eqnarray}
Since the covering $( {O_{j, \,\delta } } )_{j = 1,...,J} $ of
$\mathbb{S}^n$ depends on $\delta=\delta(\varepsilon)$, we can
choose $\delta$ such that
\begin{eqnarray}\label{E6}
\left( A_0(k)+\frac{\varepsilon_0}{2} \right)(1+\max
\varepsilon_j) \leq A_0(k)+\varepsilon_0
\end{eqnarray}
and since $\varepsilon_0$ is an arbitrary small positive real we
can  choose it such that
\begin{eqnarray}\label{E7}
\left(A_0(k)+\varepsilon_0 \right)\omega_{m-1}^{-\frac{2}{k}}\leq
A_0(k)\omega_{m-1}^{-\frac{2}{k}}+\varepsilon
\end{eqnarray}
By (\ref{E5}), because of (\ref{E6}) and (\ref{E7}), arises
\begin{eqnarray*}
\left(\int_{\mathbb{S}^n} {(\eta_j f)^2 ds}\right)^{1+\frac{2}{k}}
\leqslant \left( A_0(k)\omega_{m-1}^{-\frac{2}{k}}+
\varepsilon\right) \int_{\mathbb{S}^n} |\nabla(\eta_j f)|^2 ds
\left(\int_{\mathbb{S}^n} \left| \eta_j f \right| ds
\right)^{\frac{4} {{k }}}
\end{eqnarray*}
Last inequality means that Lemma \ref{L3.1} holds in
$\mathbb{S}^n$, and thus by Lemma \ref{L3.2} we obtain
\begin{eqnarray}\label{E8}
\left(\int_{\mathbb{S}^n} {f^2 ds}\right)^{1+\frac{2}{k}}
&\leqslant &\left(A_0(k)\omega_{m-1}^{-\frac{2}{k}}+
\varepsilon\right)\left( {\int_{\mathbb{S}^n} {\left| {\nabla f}
\right|^2 ds } + B_\varepsilon\int_{\mathbb{S}^n} {f^2 ds} }
\right) \nonumber\\
&&\times\left({\int_{\mathbb{S}^n} {\left| f \right|} ds }
\right)^{\frac{4} {{k }}}
\end{eqnarray}

We have now  to prove that the constant
$A_0(k)\omega_{m-1}^{-\frac{2}{k}}$ is the best for this
inequality. The proof of this part proceeds by contradiction,
based on inequality (\ref{E8}). We assume that, for any
$\alpha>0$, there exists $f \in C_0^\infty \left(\mathbb{S}^n
\right)$ such that
\begin{eqnarray}\label{E9}
\frac{\left( {\int_{\mathbb{S}^n} {\left| \nabla f \right|^2 ds }
+\alpha\int_{\mathbb{S}^n} {f^2 ds} } \right)
\left({\int_{\mathbb{S}^n} {\left| f\right|} ds }
\right)^{\frac{4} {{k }}}}{\left(\int_{\mathbb{S}^n} {f^2
ds}\right)^{1+\frac{2}{k}}} <
\frac{1}{A_0(k)\omega_{m-1}^{-\frac{2}{k}}}
\end{eqnarray}
Using formulas $(1.1)$ and $(1.2)$ in \cite{Ili} we can write
\begin{eqnarray}\label{E10}
\int_{\mathbb{S}^n} \left| \nabla f \right|^2 ds=
\omega_{k-1}\omega_{m-1}\int_0^{+\infty} \left| {f'}(t)(1+t^2)
\right|^2(1+t^2)^{-\frac{k+m}{2}} t^{k-1}dt
\end{eqnarray}
\begin{eqnarray}\label{E12}
\int_{\mathbb{S}^n} \left| f \right| ds=
\omega_{k-1}\omega_{m-1}\int_0^{+\infty} \left|
f(t)\right|(1+t^2)^{-\frac{k+m}{2}} t^{k-1}dt\;\;\;\;\;\;\;
\end{eqnarray}
and
\begin{eqnarray}\label{E11}
\int_{\mathbb{S}^n}  f^2ds=
\omega_{k-1}\omega_{m-1}\int_0^{+\infty}f^2(t)(1+t^2)^{-\frac{k+m}{2}}
t^{k-1}dt\;\;\;\;\;\;\;\;\;
\end{eqnarray}
For any $\lambda>0$, define $f_\lambda$ by $f_\lambda(t)=f(\lambda
t)$. So, for any $\lambda>0$, $f_\lambda$ has compact support and
by (\ref{E10}), (\ref{E11}) and (\ref{E12}) we obtain,
respectively
\begin{eqnarray}\label{E13}
\int_{\mathbb{S}^n} \!\left| \nabla f \right|^2 ds=\lambda^{2-k}
\omega_{k-1}\omega_{m-1}\!\int_0^{+\infty}\! \left|
{f'}(t)\left(1+\frac{t^2}{\lambda^2}\!\right)\!
\right|^2\!\!\left(1+\frac{t^2}{\lambda^2}\right)^{-\frac{k+m}{2}}\!\!
t^{k-1}dt\nonumber\\
\end{eqnarray}
\begin{eqnarray}\label{E15}
\int_{\mathbb{S}^n} \left| f \right| ds=\lambda^{-k}
\omega_{k-1}\omega_{m-1}\int_0^{+\infty} \left|
f(t)\right|\left(1+\frac{t^2}{\lambda^2}\right)^{-\frac{k+m}{2}}
t^{k-1}dt\;\;\;\;\;\;\;
\end{eqnarray}
and
\begin{eqnarray}\label{E14}
\int_{\mathbb{S}^n}  f^2ds=\lambda^{-k}
\omega_{k-1}\omega_{m-1}\int_0^{+\infty}f^2(t)\left(1+\frac{t^2}{\lambda^2}\right)^{-\frac{k+m}{2}}
t^{k-1}dt\;\;\;\;\;\;\;\;\;
\end{eqnarray}
For $\lambda\rightarrow +\infty$, by Lebesque's Theorem and
because of (\ref{E13}), (\ref{E14}) and (\ref{E14}), inequality
(\ref{E9}) yields
\begin{eqnarray*}
\frac{\left( { \omega_{k-1}\omega_{m-1}\int_0^{+\infty} \left|
{f'}(t) \right|^2 t^{k-1}dt} \right) \left(
\omega_{k-1}\omega_{m-1}\int_0^{+\infty} \left| f(t)\right|
t^{k-1}dt \right)^{\frac{4} {{k }}}}{\left(
\omega_{k-1}\omega_{m-1}\int_0^{+\infty}f^2(t)
t^{k-1}dt\right)^{1+\frac{2}{k}}} <
\frac{1}{A_0(k)\omega_{m-1}^{-\frac{2}{k}}}\nonumber\\
\end{eqnarray*}
or
\begin{eqnarray}\label{E16}
\frac{\left( { \omega_{k-1}\int_0^{+\infty} \left| {f'}(t)
\right|^2 t^{k-1}dt } \right) \left( \omega_{k-1}\int_0^{+\infty}
\left| f(t)\right| t^{k-1}dt \right)^{\frac{4} {{k }}}}{\left(
\omega_{k-1}\int_0^{+\infty}f^2(t)
t^{k-1}dt\right)^{1+\frac{2}{k}}}<
\frac{1}{A_0(k)}\nonumber\\
\end{eqnarray}
Since for a radial function $f$ hold
$$
f(x)=f(|x|)=f(t)\quad \mathrm{and}\quad|\nabla f(x)|=|f'(t)|
$$
we have
\begin{eqnarray}\label{E17}
\int_{\mathbb{R}^k} \left| \nabla f(x) \right|^2 dx=
\omega_{k-1}\int_0^{+\infty} \left| {f'}(t) \right|^2 t^{k-1}dt
\end{eqnarray}
Moreover, the following equalities hold
\begin{eqnarray}\label{E18}
\int_{\mathbb{R}^k} \left| {f}(x) \right| dx=
\omega_{k-1}\int_0^{+\infty} \left| {f}(t) \right| t^{k-1}dt
\end{eqnarray}
and
 \begin{eqnarray}\label{E19}
\int_{\mathbb{R}^k}  {f^2}(x) dx=
 \omega_{k-1}\int_0^{+\infty}  {f^2}(t)
 t^{k-1}dt
\end{eqnarray}
By (\ref{E16}), because of (\ref{E17}), (\ref{E18}) and
(\ref{E19}), arises
\begin{eqnarray*}
\frac{ \int_{\mathbb{R}^k} \left| \nabla f(x) \right|^2 dx
\left({\int_{\mathbb{R}^k} \left| {f}(x) \right| dx }
\right)^{\frac{4} {{k }}}}{\left(\int_{\mathbb{R}^k}  {f^2}(x)
dx\right)^{1+\frac{2}{k}}}< \frac{1}{A_0(k)}
\end{eqnarray*}
Last inequality is false (see \cite{Car-Los}) and the theorem is
proved.\mbox{ }\hfill $\Box$
 \bigbreak \noindent{\bf{Proof of Theorem \ref{T2.6}}.} The
proof of this Theorem is similar to Theorem \ref{T2.2}}.\mbox{
}\hfill $\Box$

\end{document}